\documentclass[10pt]{article}

\usepackage{amsmath}
\usepackage{amsfonts}
\usepackage{amssymb}
\usepackage{framed}
\usepackage{mathtools}
\usepackage{enumerate}
\usepackage{mathrsfs}
\usepackage[hidelinks]{hyperref}
\usepackage{enumitem}
\usepackage{mathrsfs}
\usepackage{scrextend}
\usepackage{amsthm}
\usepackage{bbm}
\usepackage{thmtools}
\usepackage{thm-restate}
\usepackage{mathabx}
\usepackage[margin=1.5in]{geometry}
\usepackage{rotating}

\usepackage{xcolor}
\hypersetup{
    colorlinks,
    linkcolor={red!50!black},
    citecolor={blue!50!black},
    urlcolor={blue!80!black}
}

\usepackage{tikz-cd}

\theoremstyle{plain}
\newtheorem{thm}{Theorem}[section]
\newtheorem{lem}[thm]{Lemma}
\newtheorem{prop}[thm]{Proposition}

\theoremstyle{definition}
\newtheorem{defn}[thm]{Definition}

\theoremstyle{remark}
\newtheorem*{rem}{Remark}

\newtheorem*{notn}{Notation}

\tikzset{
  symbol/.style={
    draw=none,
    every to/.append style={
      edge node={node [sloped, allow upside down, auto=false]{$#1$}}}
  }
}

\newcommand{\Spec}{\textrm{Spec} \hspace{0.15em} }

\newcommand\restr[2]{{
	\left.\kern-\nulldelimiterspace
	#1
	\vphantom{\big|}
	\right|_{#2}
	}}
\newcommand{\an}[1]{#1^{\textrm{an}}}

\newcommand{\Sch}{\textrm{Sch}}

\newcommand{\Hom}{\textrm{Hom}}

\newcommand{\ch}[1]{\widecheck{{#1}}}

\newcommand{\NL}{\textrm{NL}}

\newcommand{\GL}{\textrm{GL}}

\usepackage{amsmath,calligra,mathrsfs}

\usepackage[cal=boondoxo]{mathalfa}

\DeclareMathOperator{\sheafhom}{\mathcal{H \kern -1pt o \kern -2pt m}}
\DeclareMathOperator{\sheafend}{\mathcal{E \kern -1pt n \kern -2pt d}}
\DeclareMathOperator{\sheafaut}{\mathcal{A \kern -1pt u \kern -2pt t}}

\title{Unlikely Intersections with Bruhat Strata}
\author{David Urbanik}

\begin{document}

\maketitle

\begin{abstract}
Let $\mathcal{A}_{g}$ be the moduli space of $g$-dimensional principally polarized abelian varieties over $\mathbb{Z}$, and let $\mathcal{T} \subset \mathcal{A}_{g}$ be a closed locus, also defined over $\mathbb{Z}$. Motivated by unlikely intersection conjectures, we study the intersection of $\mathcal{T}_{\mathbb{F}_{p}}$ with the Bruhat strata in $\mathcal{A}_{g,\mathbb{F}_{p}}$ as $p$-varies; these are strata characterized by the existence of certain subgroup schemes inside the $p$-torsion of the fibres. We find that, away from a finite set of primes, positive-dimensional ``unlikely'' intersections of $\mathcal{T}_{\mathbb{F}_{p}}$ with such strata are all accounted for by intersections of $\mathcal{T}$ with special loci inside $\mathcal{A}_{g}$. This result generalizes to all abelian-type Shimura varieties, and variations of Hodge structures equipped with certain motivic data. It moreover gives another example of how functional transcendence principles in characteristic zero can be used to study unlikely intersections in positive characteristic, building on recent work by the author. 
\end{abstract}

\tableofcontents

\section{Introduction}
\label{intro}

The past decade of work in number theory and algebraic geometry has seen a resurgence of interest in \emph{functional transcendence principles}. The general setup is something like the following. One considers a complex algebraic variety $S$, a complex analytic variety $\widetilde{S}$ embedded as a locally analytic subspace of a complex algebraic variety $\ch{D}$, and an analytic uniformization map $\pi : \widetilde{S} \to S$. The embedding $\widetilde{S} \subset \ch{D}$ defines a notion of an algebraic subvariety of $\widetilde{S}$, which is simply an intersection $V \cap \widetilde{S}$, where $V \subset \ch{D}$ is an algebraic subvariety.  A functional transcendence principle then describes which algebraic subvarieties of $\widetilde{S}$ are of the form $\pi^{-1}(Z)$, where $Z \subset S$ is an algebraic subvariety, and moreover may characterize when such an inverse image $\pi^{-1}(Z)$ intersects ``atypically'' with subvarieties $V \subset \ch{D}$. 

Such principles are of interest in number theory and algebraic geometry because often such a setup arises in a situation where $S$ is a variety of algebraic moduli for algebro-geometric objects, and $\widetilde{S}$ parametrizes cohomological data associated to those algebro-geometric objects. In the most typical situation, one has a a smooth projective family $f : X \to S$ of complex algebraic varieties, and $\widetilde{S}$ describes those Hodge flags appearing in the cohomology of the fibres of $f$ as a subset of a flag variety $\ch{D}$. As an alternative perspective on the same thing when the cohomological data is of Shimura type, one can also consider Shimura data $(G, D)$ and Shimura varieties $S = \Gamma \backslash D$, with $\Gamma \subset G(\mathbb{Q})$ a level subgroup, and then the uniformizations $\pi : D \to S$ also admit a functional transcendence theory. In both settings this theory is now fully understood; here the relevant work is the Ax-Schanuel Theorems of \cite{AXSCHAN} and \cite{axschanshimura}, respectively. 

In the situation where the family $f : X \to S$ admits a model over the ring of integers $\mathcal{O}_{K}$ of a number field $K$, one would like to have some kind of functional transcendence principle which relates algebro-geometric subvarieties of $S_{\kappa(\mathfrak{p})}$, with $\mathfrak{p} \in \Spec \mathcal{O}_{K}$ a prime, and algebro-geometric relations on cohomological objects appearing in the fibres of the reduction $f_{\kappa(\mathfrak{p})} : X_{\kappa(\mathfrak{p})} \to S_{\kappa(\mathfrak{p})}$. But as there is no analytic geometry in positive characteristic to make use of, there is no clear way to formulate one.

In recent work \cite{algcycleloci}, we described a way to achieve functional transcendence-like results in positive characteristic, at least away from a finite subset of $\Spec \mathcal{O}_{K}$. First, one can give the variety $\ch{D}$ an integral structure, and interpret algebraic relations on cohomological objects appearing in the fibres of $f_{\kappa(\mathfrak{p})}$ in terms of algebraic subvarieties of $\ch{D}_{\kappa(\mathfrak{p})}$. One then speaks of ``infinitesimal'' correspondences between algebraic subvarieties of $\ch{D}_{\kappa(\mathfrak{p})}$ and algebraic subvarieties of $S_{\kappa(\mathfrak{p})}$, at least up to some infinitesimal ``order'' $r(\mathfrak{p})$ less than the characteristic of $\kappa(\mathfrak{p})$. If one can accumulate infinitely many of these correspondences for an infinite set of primes $\{ \mathfrak{p}_{i} \}_{i = 1}^{\infty}$ and with $\{ r(\mathfrak{p}_{i}) \}_{i = 1}^{\infty}$ unbounded, then one can often produce such a correspondence in characteristic zero, and a functional transcendence principle in characteristic zero can then be used to understand the phenomenon in question for all but finitely many of the primes in $\{ \mathfrak{p}_{i} \}_{i = 1}^{\infty}$.

In \cite{algcycleloci}, this idea was used to understand families of ordinary algebraic cycles appearing in the fibres of the families $f_{\kappa(\mathfrak{p})}$ as $\mathfrak{p}$ varied. In this paper we use these ideas, and the results of \cite{algcycleloci}, to understand a purely positive characteristic phenomenon: Bruhat strata.

\subsection{The Shimura Setting}

The main results are both easier to state and understand in the setting of Shimura varieties, so we begin with this case first. In fact, we start with the special case where $S = \mathcal{A}_{g,n}$ is the moduli space of principally polarized abelian varieties with level $n$ structure. Then we have a universal family $f : X \to S$ defined over $\mathbb{Z}$. 

\begin{defn}
The $k$'th \emph{Bruhat stratum} at the prime $p \in \mathbb{Z}$ is the locally closed subvariety $S_{k,p} \subset S_{\mathbb{F}_{p}}$ satisfying
\[ S_{k,p}(\overline{\mathbb{F}_{p}}) = \{ s \in S(\overline{\mathbb{F}_{p}}) : \dim_{\overline{\mathbb{F}_{p}}} \Hom(\alpha_{p}, X_{s}) = k \} ,\]
where $\alpha_{p} = \ker \left[ \textrm{Frob} : \mathbb{G}_{a} \to \mathbb{G}_{a} \right]$. 
\end{defn}

The subvarieties $\{ S_{k,p} \}_{k = 0}^{g}$ give an algebraic stratification of $\mathcal{A}_{g,n,\mathbb{F}_{p}}$, and the dimension of the stratum $S_{k,p}$ is $(g(g+1) - k(k+1))/2$, as shown in \cite[Ex. 4.3]{bruhatPELtype}. Supposing now we have a closed irreducible integral locus $\mathcal{T} \subset \mathcal{A}_{g,n}$, one in some sense ``expects'' that the intersection $\mathcal{T}_{\mathbb{F}_{p}} \cap S_{k,p}$ will have dimension at most $\dim \mathcal{T}_{\mathbb{F}_{p}} - k(k+1)/2$, although of course one can easily construct $\mathcal{T}$ for which the intersection is larger for finitely many $p$. Our first result shows that, should the intersection be larger for an infinite set of primes $\mathcal{P} \subset \mathbb{Z}$, then this is only on account of exceptional intersections between $\mathcal{T}$ and so-called ``special'' loci in $\mathcal{A}_{g,n}$: 

\begin{notn}
Given an algebraic variety $V$, we write $V^{> e}$ for the union of components of $V$ of dimension greater than $e$. 
\end{notn}

\begin{thm}
\label{abvarcase}
Fix a non-negative integer $e \geq \dim \mathcal{T}_{\mathbb{C}} - k(k+1)/2 .$ Then if $\mathcal{Z}$ is the Zariski closure of $\bigcup_{p \in \Spec \mathbb{Z}} [\mathcal{T}_{\mathbb{F}_{p}} \cap S_{k,p}]^{> e}$ inside $\mathcal{A}_{g,n}$, each component of $\mathcal{Z}_{\mathbb{C}}$ lies in a proper special subvariety of $\mathcal{A}_{g,n,\mathbb{C}}$.
\end{thm}

Let us explain the term ``special subvariety''. We will give a definition in the context of Shimura varieties; for a Hodge-theoretic definition see \cite[Def. 1.2]{2021InMat.225..857K}. Suppose that $X = \Gamma \backslash D$ is a Shimura variety associated to a Shimura datum $(G, D)$, where $\Gamma \subset G(\mathbb{Q})$ is an appropriate discrete subgroup. Given a Shimura subdatum $(H, E) \subset (G, D)$ and an associated Shimura variety $Y = \Gamma' \backslash E$ with $\Gamma' \subset \Gamma$, one obtains a map of Shimura varieties $t : Y \to X$. A special subvariety of $X$ is a complex algebraic subvariety obtained as the irreducible component of a Hecke translate of the image of such a $t$. In the case where $X = \mathcal{A}_{g,n,\mathbb{C}}$, these varieties conjecturally parametrize loci in $X$ where the fibres of the universal family over $X$ acquire additional algebraic cycles. 

\vspace{0.5em}

\autoref{abvarcase} is also easy to generalize to the case of abelian type Shimura varieties and their integral canonical models, at least once the appropriate notion of Bruhat stratum has been defined. In essence, such Shimura varieties $X$ and their models $\mathcal{X}$ admit correspondences $t : \mathcal{X} \dashrightarrow \mathcal{A}_{g,n}$, for appropriate $g$ and $n$, and the Bruhat strata of $\mathcal{X}_{\kappa(\mathfrak{p})}$, with $\mathfrak{p}$ some prime over which $\mathcal{X}$ is defined, are defined in terms of abelian variety data pulled back along $t_{\kappa(\mathfrak{p})}$. These strata detect both the presence of copies of the subgroup scheme $\alpha_{p}$ in the cohomology of the associated abelian variety but also take into account the extra ``$G$-structure'' present, where $(G, D)$ is the Shimura datum defining $X$; we give more details in both the next section and the proof that appears in \S\ref{shimappsec}. With this in mind, let $\mathcal{X}$ be an integral canonical model over a number ring $R = \mathcal{O}_{K}[1/N] \subset \overline{\mathbb{Q}}$, with $K$ a number field, of an abelian type Shimura variety $X = \Gamma \backslash D$ associated to $(G, D)$. We assume $G$ is a $\mathbb{Q}$-simple group, and $\Gamma \subset G(\mathbb{Q})$ a neat arithmetic subgroup. Let $\mathcal{X}_{\alpha,\mathfrak{p}} \subset \mathcal{X}_{\kappa(\mathfrak{p})}$ be an associated Bruhat stratum for all $\mathfrak{p} \in \Spec R$, and let $\mathcal{T} \subset \mathcal{X}$ be a closed algebraic locus. We then have the following:

\begin{thm}
\label{genshimcase}
Fix a non-negative integer $e \geq 0$ and let $\mathcal{P} \subset \Spec R$ be an infinite set of primes $\mathfrak{p}$ such that
\[\dim (\mathcal{T}_{\kappa(\mathfrak{p})} \cap \mathcal{X}_{\alpha,\mathfrak{p}}) > e \geq \dim \mathcal{T}_{\mathbb{C}} + \dim \mathcal{X}_{\alpha,\mathfrak{p}} - \dim X \] 
for all $\mathfrak{p} \in \mathcal{P}$. Then if $\mathcal{Z}$ is the Zariski closure of $\bigcup_{\mathfrak{p} \in \mathcal{P}} [\mathcal{T}_{\kappa(\mathfrak{p})} \cap \mathcal{X}_{\alpha,\mathfrak{p}}]^{> e}$ inside $\mathcal{X}$, each component of $\mathcal{Z}_{\mathbb{C}}$ lies in a proper special subvariety of $X$.
\end{thm}

\subsection{The Hodge-Theoretic Setting}
\label{hodgetheoreticintro}

To state our main technical result, we adopt a more general setup. Suppose that $f : X \to S$ is a smooth projective family with geometrically connected fibres, all defined over a subring $R \subset \mathbb{C}$, with $S$ a smooth base, and such that the fraction field of $R$ is a number field. We consider the following cohomological objects associated to $f$ in degree $w$:

\begin{itemize}
\item[-] The local system $\mathbb{V} = R^{w} \an{f}_{*} \mathbb{Z} / \textrm{tor}.$ modulo torsion.
\item[-] The algebraic de Rham bundle $\mathcal{H} = R^{w} f_{*} \Omega^{\bullet}_{X/S}$, with $\Omega^{\bullet}_{X/S}$ the complex of relative differentials, and equipped with the canonical isomorphism $\an{\mathcal{H}} \simeq \mathbb{V} \otimes \mathcal{O}_{\an{S}}$.
\item[-] The Hodge filtration $F^{\bullet}$ on $\mathcal{H}$, which is defined as 
\[ F^{i} \mathcal{H} = \textrm{image} \left[ R^{w} f_{*} \Omega^{\bullet \geq i}_{X/S} \to R^{w} f_{*} \Omega^{\bullet}_{X/S} \right] . \]
\item[-] The Gauss-Manin connection $\nabla : \mathcal{H} \to \Omega^{1}_{X/S} \otimes \mathcal{H}$, constructed as in \cite{katz1968}.
\item[-] For each prime $\mathfrak{p} \in \Spec R$, the \emph{conjugate} filtration $F^{\bullet}_{c}$ of the bundle $\mathcal{H}_{\kappa(\mathfrak{p})}$. To construct it one considers the filtration $\tau_{\bullet}$ on $\Omega^{\bullet}_{X/S}$ given by 
\[ \tau_{p}(\Omega^{\bullet}_{X/S}) = 0 \to \Omega^{1}_{X/S} \to \cdots \to \Omega^{p-1}_{X/S} \to \textrm{ker}\ d \to 0 \to \cdots \to 0 .  \]
This filtration induces a spectral sequence
\[ E^{p,q}_{2} = R^{p} f_{*} ( h^{q}(\Omega^{\bullet}_{X/S})) \implies R^{p+q} f_{*}(\Omega^{\bullet}_{X/S}) . \]
Base changing to $\kappa(\mathfrak{p})$ the resulting filtration on $\mathcal{H}_{\kappa(\mathfrak{p})}$ we denote by $F^{w-\bullet}_{c}$, suppressing the dependence on $\mathfrak{p}$ in the notation. We note that unlike much of the literature, our conjugate filtration is a \emph{descending} filtration. 
\end{itemize}

As a mechanism for understanding Bruhat loci in $S_{\kappa(\mathfrak{p})}$ for varying primes $\mathfrak{p} \in \Spec R$, we will construct, for closed points $s \in S$, algebraic stratifications on certain flag subvarieties $\ch{D}_{s} \subset \ch{L}_{s}$ determined by cohomological data associated to $f$; here we have $\ch{L}_{s} = \ch{L}[\mathcal{H}_{s}]$, where we define $\ch{L}[V]$ for any $\kappa$-vector space $V$ to be the $\kappa$-variety whose functor of points on $\kappa$-algebras $A$ is given by
\[ \ch{L}[V](A) = \left\{ \textrm{graded locally-free filtrations }L^{\bullet}\textrm{ of }A \otimes_{\kappa} V\textrm{ with Hodge numbers }h^{0}, \hdots, h^{w}   \right\} . \]
To describe the varieties $\ch{D}_{s}$ and other associated data, we need to introduce some Hodge-theoretic data associated to $\mathbb{V}$, for which we will need some notation. For any free $A$-module $V$ we set $\mathcal{T}(V) = \bigoplus_{a, b \geq 0} V^{\otimes a} \otimes (V^{*})^{\otimes b}$ to be the space of all tensors associated to $V$, where duality is taken with respect to $A$. We then define:

\begin{defn}
Suppose that $V$ is a free $A$-module, and $L^{\bullet}$ is a filtration on $V$ of weight $w$. Then
\begin{itemize}
\item[-] if $w$ is even, we set $L^{\textrm{mid}} = L^{w/2}$; and
\item[-] if $w$ is odd, we set $L^{\textrm{mid}} = 0$.
\end{itemize}
We will also use $L^{\textrm{mid}}$ to denote $\bigoplus_{a, b \geq 0} L^{\textrm{mid}} (V^{\otimes a} \otimes (V^{*})^{\otimes b})$, where we also write $L^{\bullet}$ for the functorially-induced filtrations on the tensor powers.
\end{defn}

\begin{defn}
The algebraic monodromy group $\mathbf{H}_{S}$ of $\mathbb{V}$ is defined to be the identity component of the Zariski closure in $\GL(\mathbb{V}_{s})$, for some $s \in S(\mathbb{C})$, of the monodromy representation
\[ \pi_{1}(S(\mathbb{C}), s) \to \GL(\mathbb{V}_{s}) \]
induced by $\mathbb{V}$. We regard this as an abstract algebraic group $\mathbf{H}_{S}$, with natural realizations $\mathbf{H}_{S,s}$ in each fibre of $\mathbb{V}$. 

More generally, if $Z \subset S_{\mathbb{C}}$ is any irreducible complex subvariety containing $s$, we write $\mathbf{H}_{Z}$ for the identity component of the Zariski closure in $\GL(\mathbb{V}_{s})$ of the image of $\pi_{1}(Z^{\textrm{nor}}(\mathbb{C}), s)$, where $Z^{\textrm{nor}} \to Z$ is the normalization. 
\end{defn}

Now recall that, from the Theorem of the Fixed Part \cite{Steenbrink1985}, the space of global sections $\mathbb{V}_{\mathbb{Q}}(S(\mathbb{C}))$ carries a natural Hodge structure which makes the maps $\mathbb{V}_{\mathbb{Q}}(S(\mathbb{C})) \to \mathbb{V}_{\mathbb{Q}, s}$ into morphism of Hodge structures. Denote by $\textrm{Hg}(\mathbb{V}_{\mathbb{Q}})$ the $\mathbb{Q}$-vector space of global Hodge tensors: i.e., elements of $\bigoplus_{a, b \geq 0} \mathbb{V}^{a,b}_{\mathbb{Q}}(S(\mathbb{C}))$, with $\mathbb{V}^{a,b}_{\mathbb{Q}} = \mathbb{V}^{\otimes a}_{\mathbb{Q}} \otimes (\mathbb{V}^{*}_{\mathbb{Q}})^{\otimes b}$, which are Hodge for this global Hodge structure. Fix a set of elements $\{ \tau_{i} \}_{i = 1}^{b} \subset \textrm{Hg}(\mathbb{V}_{\mathbb{Q}})$ such that the stabilizer subgroup $\mathbf{G}_{\tau}$ of the set $\{ \tau_{i} \}_{i = 1}^{b}$ is exactly the subspace $\textrm{Hg}(\mathbb{V}_{\mathbb{Q}})$; we regard $\mathbf{G}_{\tau}$ as an abstract group with a realization in each fibre of $\mathbb{V}$. After scaling, we can and do assume that each $\tau_{i}$ is in fact defined over $\mathbb{Z}$, and has non-zero reduction modulo each prime of $\Spec R$. The data $(\mathbf{G}_{\tau}, \mathbf{H}_{S}, F^{\bullet})$ induces, in a natural way, an abstract flag variety $\ch{D}$ defined over a number field, which we review in \S\ref{algfibbundsec}. This flag variety carries a natural open complex submanifold $D \subset \ch{D}$ and admits an analytic period map $\varphi : S_{\mathbb{C}} \to \Gamma \backslash D$ for an appropriate group $\Gamma$. It moreover has a realization in each fibre above points $s \in S(\mathbb{C})$, and these realizations we denote by $\ch{D}_{s}$. We note that in the case where $\mathbf{H}_{S}$ is the derived subgroup of $\mathbf{G}_{\tau}$ we may take $\ch{D}_{s}$ to be the orbit $\mathbf{G}_{\tau}(\mathbb{C}) \cdot F^{\bullet}_{s}$ inside $\ch{L}_{s}$. 

If one believes the Hodge conjecture, the elements of $\textrm{Hg}(\mathbb{V}_{\mathbb{Q}})$ are induced by ($\mathbb{Q}$-linear combinations of) families of algebraic cycles over $S_{\mathbb{C}}$, which one can moreover show admit models over $S_{\overline{\mathbb{Q}}}$. This means that, at least after replacing $R \subset \mathbb{C}$ with a finite integral extension and passing to an open subscheme of $\Spec R$, each $\tau_{i}$ should correspond to a global algebraic section $\nu_{i}$ of $\bigoplus_{a, b \geq 0} \mathcal{H}^{a,b}$, where $\mathcal{H}^{a,b} = \mathcal{H}^{\otimes a} \otimes (\mathcal{H}^{*})^{\otimes b}$, with the following properties (the second being satisfied for all finite primes $\mathfrak{p} \in \Spec R$):

\begin{itemize}[leftmargin=0.75in]
\item[(Flat)] We have $\nabla \nu_{i} = 0$ (see \cite[\S6.1]{algcycleloci}).
\item[($\mathfrak{p}$-Hodge)] The restriction of $\nu$ to $\mathcal{H}_{\kappa(\mathfrak{p})}$ lies inside $F^{\textrm{mid}} \cap F^{\textrm{mid}}_{c}$ (see \cite[\S6.3]{algcycleloci}).  
\end{itemize}
We will refer to these properties in the statement of \autoref{bruhatcritthm}, below.

\vspace{0.5em}

To state our technical results, we now need to explain what we mean by Bruhat loci in the reductions $S_{\kappa(\mathfrak{p})}$. Given a point $x \in S_{\kappa(\mathfrak{p})}$, we write $(\mathcal{H}_{x}, F^{\bullet}_{x}, F^{\bullet}_{c,x}, \{ \nu_{i,x} \}_{i=1}^{b} )$ for the data in the fibre above $x$. We can classify quadruples of this type into isomorphism classes, where $(H, L^{\bullet}, L^{\bullet}_{c}, \{ \mu_{i} \}_{i=1}^{b})$ is isomorphic to $(H', L'^{\bullet}, L'^{\bullet}_{c}, \{ \mu'_{i} \}_{i=1}^{b})$ if there is a linear map $H \to H'$ preserving the filtrations and tensors. If we consider such quadruples $(\mathcal{H}_{x}, F^{\bullet}_{x}, F^{\bullet}_{c,x}, \{ \nu_{i,x} \}_{i=1}^{b} )$ arising from the data above then we will see that there are only finitely many isomorphism classes, indexed by a finite set we call $\mathcal{I}$, and $\mathcal{I}$ can be chosen to be independent of both $\mathfrak{p}$ and $x$. One then obtains an algebraically-constructible \emph{partition} (not in general a stratification) 
\[ \bigsqcup_{\alpha \in \mathcal{I}} S_{\kappa(\mathfrak{p}), \alpha} \hookrightarrow S_{\kappa(\mathfrak{p})} \]
by isomorphism type. The loci in this partition we call Bruhat loci. 

Suppose now additionally that $\ch{D}$ is defined over $R$, which can be achieved after spreading out $\ch{D}$ and replacing $R$ with an open subscheme. Then the following is proven (in a more precise form) in \S\ref{stratsec}, again after possibly enlarging $R$:

\begin{prop}
\label{Dstratprop}
Suppose that for some $\mathfrak{p} \in \Spec R$, the condition ($\mathfrak{p}$-Hodge) holds for each $\nu_{i}$, and let $\kappa = \kappa(\mathfrak{p})$. Then for each point $s \in S(\overline{\kappa})$, there exists a natural stratification $\ch{D}_{\overline{\kappa}} = \bigcup_{\alpha \in \mathcal{I}} \ch{D}_{s,\alpha}$ induced by the data $(\mathcal{H}_{s}, F^{\bullet}_{s}, F^{\bullet}_{c,s}, \{ \nu_{i,s} \}_{i=1}^{b})$, with each stratum corresponding to those flags $L^{\bullet}$ for which the triple $(\mathcal{H}_{s}, L^{\bullet}, F^{\bullet}_{c, s},, \{ \nu_{i,s} \}_{i=1}^{b})$ has isomorphism type $\alpha$. Moreover, the dimension of the stratum $\ch{D}_{s,\alpha}$ depends only on $\mathfrak{p}$ and $\alpha$, and not $s$. 
\end{prop}

\begin{defn}
\label{weakspdef}
An irreducible complex subvariety $Z \subset S_{\mathbb{C}}$ is said to be \emph{weakly special} if it is maximal among such subvarieties for its algebraic monodromy group $\mathbf{H}_{Z}$. 
\end{defn}

\begin{thm}
\label{bruhatcritthm}
Suppose that $\varphi$ is quasi-finite, and that there exists an infinite set $\{ \mathfrak{p}_{j} \}_{j = 1}^{\infty}$ of closed points in $\Spec R$, an index $\alpha$, and global sections $\{ \tau_{i} \}_{i = 1}^{b}$ and $\{ \nu_{i} \}_{i = 1}^{b}$ as above, such that
\begin{itemize}
\item[-] the sections $\{ \nu_{i} \}^{b}_{i=1}$ satisfy (Flat) and ($\mathfrak{p}_{j}$-Hodge) for each $1 \leq j < \infty$; 
\item[-] there exists a non-negative integer $e \geq 0$ such that
\[ S_{\overline{\kappa(\mathfrak{p}_{j}}),\alpha} > e \geq \dim S_{\mathbb{C}} + \dim \ch{D}_{\overline{\kappa(\mathfrak{p}_{j})}, \alpha} - \dim \ch{D}_{\mathbb{C}} \]
for all $\mathfrak{p}_{j}$. 
\end{itemize}
Then if the algebraic monodromy group $\mathbf{H}_{S}$ is $\mathbb{Q}$-simple, the union $\bigcup_{j = 1}^{\infty} [S_{\kappa(\mathfrak{p}_{j}), \alpha}]^{> e}$ is not Zariski dense in $S$, and its Zariski closure $\mathcal{Z}$ satisfies the property that the components of $\mathcal{Z}_{\mathbb{C}}$ each lie inside a proper weakly subvariety of $S_{\mathbb{C}}$. 
\end{thm}

\begin{rem}
Because the algebraic monodromy group $\mathbf{H}_{S}$ is $\mathbb{Q}$-simple, it follows from \cite[Lem. 2.5]{fieldsofdef} that the weakly special subvarieties referred to in \autoref{bruhatcritthm} can in fact be taken to be special subvarieties in the sense of \cite[Def. 1.2]{2021InMat.225..857K}, which agrees with the notion of special in the Shimura case. 
\end{rem}

\begin{rem}
The assumption that $\varphi$ be quasi-finite is not essential to \autoref{bruhatcritthm}, and is imposed only to simplify the proof. One could remove this assumption by reducing to the quasi-finite case via a procedure similar to the one employed in the proof of \cite[Thm. 7.2]{algcycleloci}, but as it would add significant complexity to the argument we will instead focus on the quasi-finite case so as to not obscure the main ideas. 
\end{rem}

We note at this point that the statement of \autoref{bruhatcritthm} is equivalent to the same statement with $\Spec R$ replaced by a non-empty open subscheme (which we may again take to be the spectrum of a subring of $\mathbb{C}$), or with $R$ replaced by a finite integral extension in $\mathbb{C}$. We will use this repeatedly when proving \autoref{bruhatcritthm}, and this also justifies our assumption that the $\nu_{i}$ be defined over $R$ made above, since one can reduce a more general situation to this case by base-change and spreading out.

\section{A Criterion for Intersections with Bruhat Strata}

\subsection{Additional Setup}

\subsubsection{Algebraic Fibre Bundles over $S$}
\label{algfibbundsec}

We continue in the setting \S\ref{hodgetheoreticintro}, keeping all the notation which was introduced. Our first task will be to give an explicit construction of $\ch{D}$, and to make \autoref{Dstratprop} precise. Recall from \cite{Andre1992} that the group $\mathbf{H}_{S}$ is a $\mathbb{Q}$-normal subgroup of the reductive Mumford-Tate group $\mathbf{G}_{\tau}$. We denote by $\mathbf{N}_{S}$ its complementary normal $\mathbb{Q}$-algebraic factor, and by $\mathbf{N}_{S,s}$ its realizations in each fibre of $\mathbb{V}$. We denote by $\mathcal{P}(\mathbf{H}_{S})$ (resp. $\mathcal{P}(\mathbf{H}_{S,s})$) the $\mathbb{Q}$-variety of all parabolic subgroups of $\mathbf{H}_{S}$ (resp. $\mathbf{H}_{S,s}$). The flag $F^{\bullet}_{s}$ on $\mathbb{V}_{\mathbb{C},s}$ induces a point of $\mathcal{P}(\mathbf{H}_{S,s})$ which is the stabilizer of $F^{\bullet}_{s}$; that this stabilizer is actually parabolic can be seen directly by analyzing the adjoint Hodge structure on the Lie algebra of $\mathbf{H}_{S,s}$ induced by $F^{\bullet}_{s}$, and we will also give a different, less direct, argument below. Write $\ch{D}_{s}$ for the geometric component of $\mathcal{P}(\mathbf{H}_{S,s})$ containing this stabilizer, and let $\ch{D}$ be the corresponding component of $\mathcal{P}(\mathbf{H}_{S})$. We claim that $\ch{D}$ is well-defined (independent of $s$). This follows by analytic continuation, since on any open neighbourhood $B \subset S(\mathbb{C})$ the map $s \mapsto F^{\bullet}_{s}$ is continuous. 

The variety $\ch{D}$ is defined over a number field $E$. We now explain how to give ``algebraic de Rham'' realizations of $\mathbf{G}_{\tau}$, $\mathbf{N}_{S}$, $\mathbf{H}_{S}$ and $\ch{D}$ which are defined over $R$, at least after replacing $R$ with a finite integral extension and adjoining to $R$ finitely many denominators, and how in fact these realizations fit into smooth algebraic families over $S$. We denote by $K$ the fraction field of $R$ inside $\mathbb{C}$ in what follows. The group $\mathbf{H}_{S}$ can be realized as the identity component of the stabilizer of all global tensors of $\bigoplus_{a, b \geq 0} \mathbb{V}^{a,b}_{\mathbb{Q}}$, and this, by the Riemann-Hilbert correspondence, is identified after complexification with the space of global tensors of $\bigoplus_{a, b \geq 0} \mathcal{H}^{a,b}$ which are flat for $\nabla$. Thus, for each field-valued point $s : \Spec \kappa \to S$ we may define $\mathbf{H}_{S,s}$ to be the identity component of the stabilizer inside $\GL(\mathcal{H}_{s})$ of the image of the natural restriction 
\[ \left( \bigoplus_{a, b \geq 0} \mathcal{H}^{a,b} \right)^{\nabla = 0} \to \bigoplus_{a, b \geq 0} \mathcal{H}^{a,b}_{s} . \]
This definition agrees with the usual one at complex fibres, and also makes sense at positive characteristic fibres. To define $\mathbf{G}_{\tau,s}$ at a general field-valued point $s : \Spec \kappa \to S$, we simply take it to be the stabilizer inside $\GL(\mathcal{H}_{s})$ of $\{ \nu_{i} \}^{b}_{i=1}$. Note that we are implicitly identifying here the complex algebraic de Rham and complex Betti realizations of these groups, which is largely harmless since our setup ensures the $R$-algebraic structure on both sides of the Betti-algebraic de Rham comparison is sufficiently compatible, and context will make it clear which group we mean. 

As for the case of $\mathbf{N}_{S}$, one can define the realizations $\mathbf{N}_{S,s} \subset \GL(\mathcal{H}_{s})$ for an arbitrary $s : \Spec \kappa \to S$ to be the identity component of the stabilizer subgroup of those tensors inside $\bigoplus_{a, b \geq 0} \mathcal{H}^{a,b}_{s}$ which are either fixed by $\mathbf{G}_{\tau,s}$ or not fixed by $\mathbf{H}_{S,s}$.

Now observe that, after replacing $R$ with a further proper subring of $K$, we may assume that the groups $\mathbf{H}_{S,s}$, $\mathbf{N}_{S,s}$ and $\mathbf{G}_{\tau,s}$ fit into smooth families of connected group schemes inside the $R$-algebraic $\GL_{m,R}$-torsor $\mathcal{Aut}(\mathcal{H}) \to S$, where $\mathcal{Aut}(\mathcal{H}) \to S$ is the natural torsor above $S$ whose fibres are given by $\GL(\mathcal{H}_{s})$. Indeed, the property that a morphism be smooth is open on the target, and the required properties hold over the generic point of $\Spec R$. We may also assume that $\mathbf{H}_{S,s} \cdot \mathbf{N}_{S,s}$ is an almost-direct product decomposition of $\mathbf{G}_{\tau,s}$ at each such point $s$, again by spreading out from the generic fibre. One similarly obtains, after possibly increasing $R$, a smooth family $\ch{\mathcal{D}} \to S$ whose fibres above points $s \in S(\mathbb{C})$ are naturally identified with $\ch{D}_{s}$, and for which the fibre above any point $s : \Spec \kappa \to S$ is identified with the component in $\mathcal{P}(\mathbf{H}_{S,s})$ which contains the stabilizer of $F^{\bullet}_{s}$. 

To construct, at each positive-characteristic point $s$ of $S$, the desired stratifications of $\ch{D}_{s}$, we need to better understand the relation between the flags $F^{\bullet}_{s}$ and $F^{\bullet}_{c,s}$. Let us take $\ch{\mathcal{L}} \to S$ to be the natural family of flag schemes whose fibres are identified with $\ch{L}_{s}$. We consider the flag subbundle $\mathrm{NL}_{\tau} \subset \ch{\mathcal{L}}$ whose fibres $\mathrm{NL}_{\tau,s}$ are identified with the flag subvarieties of $\ch{L}_{s}$ consisting of flags $L^{\bullet}$ for which $\textrm{span} \{ \nu_{1,s}, \hdots, \nu_{b,s} \} \subset L^{\textrm{mid}}$. Let us elucidate the geometry of $\mathrm{NL}_{\tau,s}$ in the situation where $s$ is an element of $S(\mathbb{C})$. Recall that the variation $\mathbb{V}$ admits a natural polarization by cup product on its primitive part, which can be extended to a polarization on the whole of $\mathbb{V}$ as in \cite[Cor. 2.3.5]{CMS}. Then there is no harm in assuming that one of the $\{ \tau_{i} \}_{i = 1}^{b}$, say it is $\tau_{1}$, corresponds to exactly such a polarizing form $Q$; note that the properties (Flat) and ($\mathfrak{p}_{j}$-Hodge) of \autoref{bruhatcritthm} will then always be satisfied for $\nu_{1}$ since $\nu_{1}$ is an algebraic cycle class. Viewing, via the Betti-algebraic de Rham isomorphism, the variety $\mathrm{NL}_{\tau,s}$ inside of $\ch{L}[\mathbb{V}_{\mathbb{C},s}] \simeq \ch{L}_{s}$, it then agrees with the variety $\ch{\textrm{NL}}_{M}$ of \cite[pg.216]{GGK}, with $M = \mathbf{G}_{\tau,s}$. By \cite[VI.B.1]{GGK}, one therefore finds that the components of $\mathrm{NL}_{\tau,s}$ consist of finitely many $\mathbf{G}_{\tau,s}(\mathbb{C})$-orbits inside $\ch{L}_{s}$. Note that this in particular implies that these orbits of $\mathbf{G}_{\tau,s}$ are complete algebraic varieties.

Returning to the algebraic de Rham setting we therefore learn that the algebraic bundle $\mathrm{NL}_{\tau,\mathbb{C}} \to S_{\mathbb{C}}$ decomposes into components $\pi_{i, \mathbb{C}} : \NL_{\tau,i,\mathbb{C}} \to S_{\mathbb{C}}$ for $i = 1, \hdots, n$, where each $\pi_{i}$ is smooth and proper, and with the fibres of each $\pi_{i,\mathbb{C}}$ homogeneous spaces for $\mathbf{G}_{\tau, \mathbb{C}}$. After spreading out and replacing $R$ with a finite integral extension, we may assume the same is true about maps $\pi_{i} : \NL_{\tau,i} \to S$ defined over $R$ which give the components of $\mathrm{NL}_{\tau}$. Notice that $F^{\bullet}_{s}$ induces a section of one of the $\pi_{i}$; we assume it is $\pi_{1}$ without loss of generality. We now note that, with this setup, the stabilizer $P_{s} \subset \GL(\mathcal{H}_{s})$ of $F^{\bullet}_{s}$, with $s : \Spec \kappa \to S$ any field-valued point, induces a parabolic subgroup $P_{s} \cap \mathbf{H}_{S,s}$ of $\mathbf{H}_{S,s}$. To see this, first note that $P_{s} \cap \mathbf{G}_{\tau,s}$ is parabolic in $\mathbf{G}_{\tau,s}$, which follows because $\mathbf{G}_{\tau,s}/(P_{s} \cap \mathbf{G}_{\tau,s})$ is isomorphic to the complete variety $\NL_{\tau,1,s}$. Then because $\mathbf{G}_{\tau,s} = \mathbf{H}_{S,s} \cdot \mathbf{N}_{S,s}$ is almost-direct, one additionally has the decomposition
\[ \mathbf{G}_{\tau,s}/(P_{s} \cap \mathbf{G}_{\tau,s}) = \mathbf{H}_{S,s} / (P_{s} \cap \mathbf{H}_{S,s}) \times \mathbf{N}_{S,s} / (P_{s} \cap \mathbf{N}_{S,s}) , \]
and both factors on the right are complete varieties because the variety on the left is. From this parabolicity and factorization, we in particular obtain a canonical map $\rho: \NL_{\tau,1} \to \ch{\mathcal{D}}$ over $S$ which sends a flag $L^{\bullet} \in \NL_{\tau,1,s}(\kappa)$ to the parabolic subgroup of $\mathbf{H}_{S,s}$ which stabilizes it.

\subsubsection{Stratifications}
\label{stratsec}

As a predecessor to constructing stratifications of $\ch{D}_{s}$ for each field valued point $s : \Spec \kappa \to S$, we first begin by constructing stratifications of $\NL_{\tau,1,s}$; the eventual stratification of $\ch{D}_{s}$ we are interested in will be obtained by taking the image under $\rho$. For each such stratification we will fix a flag $L^{\bullet}_{c} \in \NL_{\tau,s}(\kappa)$ and the resulting stratification of flags $L^{\bullet} \in \NL_{\tau,1,s}(\kappa)$ will be by the isomorphism class of the quadruple $(\mathcal{H}_{s}, L^{\bullet}, L^{\bullet}_{c}, \{ \nu_{i,s} \}_{i=1}^{b})$.

\begin{lem}
\label{bruhatstratlem}
Fix an index $o$. For each fixed $L^{\bullet}_{c} \in \NL_{\tau,o,s}(\kappa)$, we have a stratification 
\[ \NL_{\tau,1,s} = \bigcup_{\alpha \in \mathcal{I}_{o}} \mathcal{S}(s, L^{\bullet}_{c}, \alpha) \subset \NL_{\tau,1,s} \]
of $\NL_{\tau,1,s}$ by the isomorphism class of the quadruple $(\mathcal{H}_{s}, L^{\bullet}, L^{\bullet}_{c}, \{ \nu_{i,s} \}_{i=1}^{b})$. The set $\mathcal{I}_{o}$ of isomorphism classes is finite and naturally identified with the double coset $W_{o} \backslash W / W_{1}$, with
\begin{itemize}
\item[-] the group $W$ being the Weyl group of $\mathbf{G}_{\tau,s}$, obtained by fixing a maximal torus $T \subset \mathbf{G}_{\tau,s}$ and taking $W = N_{\mathbf{G}_{\tau,s}}(T)/T$; and
\item[-] for each $j$ with $1 \leq j \leq n$, the subgroup $W_{j} \subset W$  is the group generated by $I_{j} \subset W$, where $I_{j}$ is the set of simple reflections lying inside some parabolic subgroup $P_{j} \subset \mathbf{G}_{\tau,s}$ containing $T$ which is the stabilizer of some point of $\NL_{\tau,j,s}$.
\end{itemize}
\end{lem}

\begin{proof}
As all the data in the quadruple $(\mathcal{H}_{s}, L^{\bullet}, L^{\bullet}_{c}, \{ \nu_{i,s} \}_{i=1}^{b})$ except for $L^{\bullet}$ is fixed, the isomorphism classes are naturally identified with the orbits of $P_{o} \subset \mathbf{G}_{\tau,s}$ inside $\NL_{\tau,1,s}$; here $P_{o}$ is the group defined by the property that it fixes $L^{\bullet}_{c}$. Since we have an identification $\NL_{\tau,1,s} = \mathbf{G}_{\tau,s} / P_{1}$, where $P_{1}$ is the stabilizer of some point of $\NL_{\tau,1,s}$, we can also identify the set of equivalence classes with the double coset space $P_{o} \backslash \mathbf{G}_{\tau,s} / P_{1}$. Since $\mathbf{G}_{\tau,s}$ acts transitively on $\NL_{\tau,1,s}$ we may choose $P_{1}$ so that $P_{o}$ and $P_{1}$ contain a common Borel subgroup $B$. The result now follows from the Bruhat decomposition \cite[14.12]{borel1991linear} of $\mathbf{G}_{\tau,s}$ for $B$.
\end{proof}

\begin{defn}
Continuing with the notation introduced in \autoref{bruhatcritthm}, for $w \in W$ we write $\ell(w)$ for the \emph{length} of $w$, which is the minimum length of a presentation of $w$ as a product of simple reflections.
\end{defn}

\begin{lem}
\label{bruhatdim}
In the notation of \autoref{bruhatcritthm}, we have 
\[ \dim \mathcal{S}(s, L^{\bullet}_{c}, \alpha) =  \textrm{max}_{w \in W_{o} w_{\alpha} W_{1}} \ell(w) , \]
where $w_{\alpha}$ is a representative of the coset corresponding to $\alpha$.
\end{lem}

\begin{proof}
From our description in the proof of \autoref{bruhatstratlem} each stratum is a finite union of orbits of $B$, so it suffices to take the maximum over the dimensions of these orbits. The required result is then implicit in the proof of \cite[21.29]{borel1991linear}, which establishes an isomorphism between Bruhat double cosets $B w B$ and a variety of dimension $\ell(w) + \dim B$. 
\end{proof}

Now suppose that $\kappa = \kappa(\mathfrak{p})$ is the residue field of a finite prime of $R$. Then using the property ($\mathfrak{p}$-Hodge), the conjugate filtration $F^{\bullet}_{c}$ associated to $\mathcal{H}_{\kappa}$ induces a section of $\pi_{o} : \NL_{\tau,o,\kappa} \to S_{\kappa}$ for some $o = o(\mathfrak{p})$ dependent on $\mathfrak{p}$, and hence a stratification 
\[ \bigsqcup_{\alpha \in \mathcal{I}_{o}} \mathcal{S}(\mathfrak{p}, \alpha) \hookrightarrow \NL_{\tau,1,\kappa} , \]
by algebraic fibre subbundles of $\pi_{1}$. Here we have $\mathcal{S}(\mathfrak{p},\alpha)_{s} = \mathcal{S}(s, F^{\bullet}_{c,s}, \alpha)$. Likewise we also obtain a stratification
\[ \bigsqcup_{\alpha \in \mathcal{I}} \rho(\mathcal{S}(\mathfrak{p}, \alpha)) \hookrightarrow \ch{\mathcal{D}} \]
by algebraic fibre subbundles of $\ch{\mathcal{D}} \to S$. We finally define the stratification of \autoref{Dstratprop} by $\ch{D}_{s,\alpha} = \ch{D}_{s} \cap \rho(\mathcal{S}(\mathfrak{p}, \alpha))$ for a point $s$ which lies over $\mathfrak{p}$. 

\subsection{Jet Theory and Differential Correspondences}

We now recall a construction developed in \cite{algcycleloci} for studying infinitesimal variations in Hodge flags using jets. For each pair of integers $d, r \geq 0$, we set $A^{d}_{r} = \mathbb{Z}[t_{1}, \hdots, t_{d}]/(t_{1}, \hdots, t_{d})^{r+1}$, and let $\mathbb{D}^{d}_{r} = \Spec A^{d}_{r}$ be the associated scheme. By a jet (or $R$-jet) of a scheme $X$ over a ring $R$, we mean a map $\mathbb{D}^{d}_{r,R} \to X$. For each $d, r \geq 0$, there are natural spaces parameterizing jets of $X$, defined as follows: 

\begin{defn}
Suppose that $S$ is an $R$-scheme, the \emph{jet space} $J^{d}_{r} S$ is defined to be the $R$-scheme representing the functor $\Sch_{R} \to \textrm{Set}$ given by
\[ T \mapsto \Hom_{R}(T \times_{R} \mathbb{D}^{d}_{r, R}, S), \hspace{1.5em} [T \to T'] \mapsto [\Hom_{R}(T' \times_{R} \mathbb{D}^{d}_{r, R}, S) \to \Hom_{R}(T \times_{R} \mathbb{D}^{d}_{r, R}, S)] , \]
where the natural map $\Hom_{R}(T' \times_{R} \mathbb{D}^{d}_{r, R}, S) \to \Hom_{R}(T \times_{R} \mathbb{D}^{d}_{r, R}, S)$ obtained by pulling back along $T \times_{R} \mathbb{D}^{d}_{r, R} \to T' \times_{R} \mathbb{D}^{d}_{r, R}$.
\end{defn}

\noindent The representability of the functor defining $J^{d}_{r} S$ reduces as in \cite[\S2.1]{periodimages} to the representability of Weil restrictions, and hence holds for instance when $S$ is quasi-projective over $R$. Note that a map $g : S \to S'$ of $R$-schemes induces a natural map $J^{d}_{r} S \to J^{d}_{r} S'$ by post-composition. Moreover, for any embedding $\mathbb{D}^{d}_{r} \hookrightarrow \mathbb{D}^{d'}_{r'}$ of disks we obtain a natural map $J^{d'}_{r'} S \to J^{d}_{r} S$ by pulling back jets; as a special case we have projections $J^{d}_{r} S \to S$ which send jets to the points at which they land.

Next, let us introduce infinitesimal period maps. We denote by $q : \GL_{m} \to \ch{L}$ the natural map which sends a basis of $A^{m}$, for $A$ a ring, to the filtration $L^{\bullet}$ for $A^{m}$ for which $L^{k} A^{m}$ is the span of the first $h^{0} + \hdots + h^{w - k}$ basis vectors. Fix an integer $r \geq 0$, and a point $s : \Spec \kappa \to S$ for a field $\kappa$ which is either of characteristic zero, or has positive characteristic $p > r$. Let $\mathfrak{m}_{s}$ be the maximal ideal of the local ring $\mathcal{O}_{S_{\kappa}, s}$ of $S_{\kappa}$ at $s$. The data $(\mathcal{H}, F^{\bullet}, \nabla)$ naturally induces a filtered modules with connection $(\mathcal{H}^{r}_{s}, F^{r, \bullet}_{s}, \nabla^{r}_{s})$ over the rings $\mathcal{O}^{r}_{S_{\kappa},s} = \mathcal{O}_{S_{\kappa},s}/\mathfrak{m}^{r+1}_{s}$ for each such $s$ and $r$. We say a frame $b^{1}, \hdots, b^{m}$ for $\mathcal{H}^{r}_{\delta}$ is flat if $\nabla^{r}_{s} b^{i} = 0$ for all $1 \leq i \leq m$. We will also allow the case $r = \infty$ when $\kappa$ has characteristic zero, in which case $\mathcal{O}^{\infty}_{S_{\kappa},s}$ is simply the formal completion of $\mathcal{O}_{S_{\kappa},s}$ at $s$.

\begin{defn}
\label{locpermapdef}
With the above setup and some $r \in \mathbb{N} \cup \{ \infty \}$, suppose that we have a filtration-compatible frame $v^{1}, \hdots, v^{m}$ for $\mathcal{H}^{r}_{s}$, and a flat frame $b^{1}, \hdots, b^{m}$ for $\mathcal{H}^{r}_{s}$. Then if $M \in \GL_{m}(\mathcal{O}^{r}_{S_{\kappa},s})$ is the change-of-basis matrix from $v^{1}, \hdots, v^{m}$ to $b^{1}, \hdots, b^{m}$ then we call the composition $\psi = q \circ M$ a \emph{local infinitesimal period map} of order $r$ at $s$. For ease of terminology we will also say $\psi$ is a ``$r$-limp''. 
\end{defn}

\begin{defn}
We say a pair $(\psi, \iota)$ consisting of an $r$-limp $\psi$ at $s$ as above and a frame $\iota : \mathcal{H}_{s} \to \kappa^{m}$ is a \emph{framed} $r$-limp if the reduction of $b^{1}, \hdots, b^{m}$ to $\mathcal{H}_{s}$ agrees with the frame induced by $\iota$. 
\end{defn}

\begin{lem}
If $r!$ is invertible in $R$ then for every point $s : \Spec \kappa \to R$ and every frame $\iota : \mathcal{H}_{s} \to \kappa^{m}$ there exists $\psi$ such that $(\psi, \iota)$ is a framed $r$-limp. 
\end{lem}

\begin{proof}
This is \cite[Lem 2.6]{algcycleloci}. 
\end{proof}

In \cite[\S2]{algcycleloci} we gave a universal construction for computing with $r$-limps associated to a triple $(\mathcal{H}, F^{\bullet}, \nabla)$ consisting of a filtered vector bundle with flat connection over a smooth $R$-scheme $S$, where $R$ is assumed to be a $\mathbb{Z}[1/r!]$-algebra. The main statement is as follows, which is a summary of \cite[Thm. 1.14]{algcycleloci} and \cite[Lem. 2.9]{algcycleloci}: 

\begin{thm}
\label{bigjetthm}
There exists a map 
\[ \eta^{d}_{r} : J^{d}_{r} S \to \GL_{m,R} \backslash J^{d}_{r} \ch{L}_{R} \]
of $R$-algebraic stacks with the following properties:
\begin{itemize}
\item[(i)] The $\GL_{m,R}$-torsor $\gamma : \mathcal{P}^{d}_{r} \to J^{d}_{r} S$ associated to $\eta^{d}_{r}$ is the base-change along the map $J^{d}_{r} S \to S$ of the frame bundle of $\mathcal{H}$.
\item[(ii)] The $\GL_{m,R}$-invariant map $\alpha : \mathcal{P}^{d}_{r} \to J^{d}_{r} \ch{L}_{R}$ associated to $\eta^{d}_{r}$ satisfies the following property. Suppose $(\iota, j) \in \mathcal{P}^{d}_{r}(\kappa)$ is a point of $\mathcal{P}^{d}_{r}$ with $\kappa$ a field, $j \in (J^{d}_{r}S)(\kappa)$ a jet lying above $s \in S(\kappa)$, and $\iota = (b^{1}_{0}, \hdots, b^{m}_{0})$ a frame of $\mathcal{H}_{s}$. Then if $(\psi, \iota)$ is a framed $r$-limp, we have $\alpha(\iota,j) = \psi \circ j$.
\item[(iii)] The construction of $\eta^{d}_{r}$ (by which we mean the pair $(\gamma, \alpha)$) is compatible with base-change along maps $R \to R'$ of $\mathbb{Z}[1/r!]$-algebras and pullback along maps $S' \to S$ of smooth quasi-projective $R$-schemes, in the sense made precise by \cite[Thm. 1.14(i)]{algcycleloci} and \cite[Thm. 1.14(iii)]{algcycleloci}. 
\end{itemize}
\end{thm}

The jet-theoretic tools will be primarily use to apply a result proven in \cite{algcycleloci} that allows us to explicitly construct a non-Zariski dense locus where atypical intersections occur. In what follows we denote by $\Phi^{d}_{r}$ the graph of the map $\alpha$ referenced in \autoref{bigjetthm}. 

\begin{defn}
For any $R$-scheme $X$ and $d, r \geq 0$, we write $J^{d}_{r,nd} X \subset J^{d}_{r} X$ for the open subspace of non-degenerate jets, defined as follows:
\begin{itemize}
\item for $r = 0$ we have $J^{d}_{0,nd} X = X$;
\item for $r = 1$ we have a natural projection $\nu : J^{d}_{1} X \to (TX)^{d}$ to the $d$-times self-product of the relative tangent bundle induced by the natural map $R[t_{1}, \hdots, t_{d}]/(t_{1}, \hdots, t_{d})^2 \to R[t_{1}]/(t_{1})^2 \times \cdots \times R[t_{d}]/(t_{d})^2$, and we define $J^{d}_{1,nd} X$ to be the fibre above the locus of $d$-tuples of linearly independent vectors above a common point of $X$;
\item for $r \geq 2$ we let $J^{d}_{r,nd} X$ be the fibre above $J^{d}_{1,nd} X$.
\end{itemize}
\end{defn}

\begin{defn}
For any $R$-scheme $X$ and $d, r \geq 0$, we write $J^{d}_{r,nc} X \subset J^{d}_{r} X$ for the open subspace of non-constant jets, defined as the complement of the natural ``zero section'' $X \to J^{d}_{r} X$ which sends a point to the constant map from the disk $\mathbb{D}^{d}_{r}$ mapping to it. 
\end{defn}

\begin{prop}[Prop. 5.8 of \cite{algcycleloci}]
\label{zariskidensityprop}
Suppose that the algebraic monodromy group $\mathbf{H}_{S}$ associated to the variation $\mathbb{V}$ underlying $(\mathcal{H}_{\mathbb{C}}, F^{\bullet}_{\mathbb{C}}, \nabla_{\mathbb{C}})$ is $\mathbb{Q}$-simple. Let $g : \mathcal{Y} \to \mathcal{M}$ be an $R$-algebraic family of geometrically irreducible subschemes of $\ch{L}_{R}$ over the Notherian base $\mathcal{M}$, with projection $\pi : \mathcal{Y} \to \ch{L}_{R}$ restricting to an embedding on fibres, and with fibres of dimension $a$. Suppose that for some integer $d$ we have
\[ \dim \ch{D}_{\mathbb{C}} - a > \dim S_{\mathbb{C}} - d . \]
Let $\mathcal{E}_{r} \subset \Phi^{d}_{r} \cap (\gamma^{-1}(J^{d}_{r,nd} S) \times J^{d}_{r,nc} \ch{L}_{R})$ be the constructible locus of points whose image inside $J^{d}_{r,nc} \ch{L}_{R}$ lies inside $J^{d}_{r} (\mathcal{Y}_{m})$ for some $m$, where the map $m : \kappa(\mathfrak{q}) \to \mathcal{M}$ is induced by a scheme-theoretic point $\mathfrak{q} \in \mathcal{M}$. Then there exists $r$ such that, after replacing $R$ with $R[1/r!]$ and applying base-change, the image in $S$ of $\mathcal{E}_{r}$ lies in a proper closed subscheme $E \subset S$. 
\end{prop}

\subsection{The Criterion}

We now relate our constructions involving stratifications to the notion of $r$-limp we have introduced. We note that if $s : \Spec \kappa \to S$ is a point and $\iota : \mathcal{H}_{s} \xrightarrow{\sim} \kappa^{m}$ is a frame, there is a canonical identification $\ch{\iota}^{-1} : \ch{L}_{\kappa} \xrightarrow{\sim} \ch{L}_{s}$ induced by $\iota$.

\begin{lem}
\label{factorthroughlem}
Fix a point $s : \Spec \kappa \to S$, where $\kappa$ is a field of positive characteristic $p$ and $s$ maps to $\mathfrak{p} \in \Spec R$, and suppose that $s$ lies inside $S_{\kappa, \alpha}$. Then if $(\psi, \iota)$ is a framed $r$-limp at $s$ (necessarily with $r < p$), the restriction of $\ch{\iota}^{-1} \circ \psi$ to $\overline{S_{\kappa, \alpha}}^{\textrm{Zar}}$ lands inside $\overline{\mathcal{S}(\mathfrak{p}, \alpha)}^{\textrm{Zar}}$. 
\end{lem}

\begin{proof}
Write $\psi = q \circ M$, where $M$ is a change-of-basis matrix between a filtration-compatible frame $v^{1}, \hdots, v^{m}$ for the pair $(\mathcal{H}^{r}_{s}, F^{r,\bullet}_{s})$ and $b^{1}, \hdots, b^{m}$ is a flat frame for $(\mathcal{H}^{r}_{s}, F^{r,\bullet}_{s}, \nabla^{r}_{s})$ which extends the frame $b^{1}_{s}, \hdots, b^{m}_{s}$ given by $\iota$. Then $\psi$ is a coordinate form of the map $\psi_{s} = \ch{\iota}^{-1} \circ \psi$, in the sense that $\psi_{s}$ admits a decomposition $\psi_{s} = q_{s} \circ M_{s}$, with $M_{s}$ the $\mathcal{O}^{r}_{S_{\kappa},s}$-linear automorphism of $\mathcal{H}^{r}_{s}$ which sends $v^{i}$ to $b^{i}$ for $1 \leq i \leq m$, and $q_{s} : \GL(\mathcal{H}_{s}) \to \ch{L}_{s}$ is the analogous quotient map. 

Now denote by $F^{r,\bullet}_{c,s}$ the restriction of $F^{\bullet}_{c}$ to $\mathcal{H}^{r}_{s}$. We claim that a $\nabla^{r}_{s}$-flat section $b$ of $\mathcal{H}^{r}_{s}$ whose fibre $b_{s}$ at $\mathcal{H}_{s}$ lies inside $F^{i}_{c,s}$ for some $i$ in fact lies inside $F^{r,i}_{c,s}$. From \cite[2.3.0]{Katz1972}, the filtration $F^{\bullet}_{c}$ satisfies the property that $\nabla F^{i}_{c} \subset \Omega^{1}_{S_{\kappa}} \otimes F^{i}_{c}$ for all $i$, so this property also holds for $F^{r,i}_{c,s}$ for all $i$. Choose an arbitrary section $v$ of $F^{r,i}_{c,s}$ whose fibre $v_{s}$ at $s$ agrees with $b_{s}$. Choose local coordinates $x_{1}, \hdots, x_{n}$ for the smooth variety $S_{\kappa}$ at $s$ such that $dx_{1}, \hdots, dx_{n}$ generate the stalk of $\Omega^{1}_{S_{\kappa}}$ at $s$. Then we may consider the operator $P : \mathcal{H}^{r}_{s} \to \mathcal{H}^{r}_{s}$ defined by (c.f. \cite[Proof of Thm. 5.1]{katznil})
\[ v' \mapsto \sum_{a} \prod_{i = 1}^{n} \left( \frac{(-x_{i})^{a_{i}}}{a_{i}!} \right) \left( \prod_{i=1}^{n} \nabla^{r}_{s} \left(\frac{\partial}{\partial x_{i}} \right)^{a_{i}} \right) (v') , \]
where the sum is over $n$-tuples $(a_{1}, \hdots, a_{n})$ of non-negative integers for which $a_{1} + \cdots + a_{n} \leq r$. Note that the fact that the expression is well-defined uses the integrability of the connection $\nabla$, as this allows us to take the derivatives in the second product in any order. By explicit computation, one sees that $P(v')$ is $\nabla^{r}_{s}$-flat for any $v'$, and the fibre of $P(v')$ at $s$ agrees with $v'_{s}$. Moreover, one similarly checks using explicit coordinates that a section of $\mathcal{H}^{r}_{s}$ is determined by the values of its $\nabla^{r}_{s}$-derivatives at $s$ up to order $r$, from which it follows that $P(v')$ depends only on $v'_{s}$. Now because $F^{i}_{c}$ is preserved by $\nabla$, one has that $P(v)$ lies inside $F^{r,i}_{c,s}$. But this means this is true for $b = P(b) = P(v)$. 

From this analysis and the flatness of the sections $\{ \nu_{i} \}_{i=1}^{b}$, it follows that under the trivialization $\beta : \mathcal{H}^{r}_{s} \xrightarrow{\sim} (\mathcal{O}^{r}_{S_{\kappa}, s})^{m}$ induced by $b^{1}, \hdots, b^{m}$ the quadruple $(\mathcal{H}^{r}_{s}, F^{r,\bullet}_{s}, F^{r,\bullet}_{c,s}, \{ \nu^{r}_{i,s} \}_{i=1}^{b})$ may be identified with $((\mathcal{O}^{r}_{S_{\kappa}, s})^{m}, \beta(F^{r,\bullet}_{s}), \mathcal{O}^{r}_{S_{\kappa}, s} \otimes L^{\bullet}_{c}, \{ 1 \otimes \mu_{i} \}_{i=1}^{b} )$, where $L^{\bullet}_{c}$ is the fixed filtration on $\kappa^{m}$ obtained as the image under $\iota$ of $F^{\bullet}_{c,s}$ and $\mu_{i}$ is the fixed tensor obtained as the image under $\iota$ of $\nu_{i,s}$. The flag $\beta(F^{r,\bullet}_{s})$ induces a map to $\ch{L}_{\kappa}$ from the formal $r$-th order neighbourhood $\mathcal{N}(r,s) \subset S_{\kappa}$ around $s$. If we restrict to $B = \mathcal{N}(r,s) \cap S_{\kappa,\alpha}$ then, by definition, the quadruple $\restr{\beta^{-1}((\mathcal{O}^{r}_{S_{\kappa}, s})^{m}, \beta(F^{r,\bullet}_{s}), \mathcal{O}^{r}_{S_{\kappa}, s} \otimes L^{\bullet}_{c}, \{ 1 \otimes \mu_{i} \}_{i=1}^{b})}{B}$ has the isomorphism type $\alpha \in \mathcal{I}_{o}$, where $o$ is an index for which $F^{\bullet}_{c,s} \in \NL_{\tau,o,s}(\kappa)$. (Note that while we have defined the stratum $S_{\kappa,\alpha}$ in the introduction by taking $\alpha$ to be the isomorphism type of its field-valued points, one could do so more generally with scheme-valued points, as is done for instance in \cite[Def. 1.13, Rem. 1.14]{bruhatstack} in much greater generality.) Concretely, this means that the infinitesimal family of parabolic subgroups of $\mathbf{G}_{\tau,s}$ induced by this quadruple lies inside a fixed double coset $P_{o} w_{\alpha} P_{1} \subset \mathbf{G}_{\tau,s}$ corresponding to $\alpha$ (in the notation of \S\ref{stratsec}), or that the map $\beta \to \ch{L}_{\kappa}$ lands inside the stratum $\ch{\iota}(\mathcal{S}(\alpha, \mathfrak{p})_{s})$. Translating back via $\ch{\iota}^{-1}$, this gives exactly what we wanted to show.
\end{proof}

We now prove our main result:

\begin{proof}[Proof of \ref{bruhatcritthm}:]
Note that the theorem statement formally reduces to the case where $a = \dim \ch{D}_{\overline{\kappa(\mathfrak{p}_{j})},\alpha}$ is assumed independent of $j$. We will construct an appropriate family $g : \mathcal{Y} \to \mathcal{M}$ as in the statement of \autoref{zariskidensityprop} (in particular with $\mathcal{M}$ a Noetherian base), and show that \autoref{zariskidensityprop} allows us to produce a proper Zariski closed subset $\mathcal{Z} \subset S$ containing the union $\bigcup_{j = 1}^{\infty} S_{\kappa(\mathfrak{p}_{j}),\alpha}$. To show the claim about $\mathcal{Z}_{\mathbb{C}}$, we will base-change the family $f : X \to S$ we started with to the components of $\mathcal{Z}$ above the generic point of $\Spec R$ and induct until the algebraic monodromy group (and hence the dimension of $\ch{D}_{\mathbb{C}}$) changes. 

Consider the set $\mathcal{V}$ of pairs $(V, s)$, where
\begin{itemize}
\item[-] $s$ is a point $s : \Spec \kappa \to S$, for some field $\kappa$, which maps to the finite prime $\mathfrak{p} \in \Spec R$ under the structure map; 
\item[-] $V$ is a closed subvariety $V \subset \ch{L}_{\kappa}$ obtained via the following process: one chooses a frame $\iota: \mathcal{H}_{s} \xrightarrow{\sim} \kappa^{m}$ and defines $V$ to be the image under the induced isomorphism $\ch{\iota} : \ch{L}_{s} \xrightarrow{\sim} \ch{L}$ of the intersection $\mathcal{S}(\mathfrak{p}, \alpha) \cap (\mathbf{H}_{S,s} \cdot F^{\bullet}_{s})$.
\end{itemize}
We wish to construct a family $g : \mathcal{Y} \to \mathcal{M}$ as above such that all such $V$ with dimension $a$ occur as fibres. In fact, it suffices to do the same for all such $V$ without this dimension condition, since one can always stratify the base $\mathcal{M}$ by the dimension of the fibres and replace $\mathcal{M}$ with a locally closed subscheme if necessary. To construct $g$, we consider the strata $\ch{\iota}(\mathcal{S}(s, L^{\bullet}_{c}, \alpha))$ ranging over all choices of $\iota$, $s$ and $L^{\bullet}_{c}$, with the notation as in \S\ref{stratsec}, and consider all possible intersections of such strata with the varieties $\ch{\iota}(\{ q \} \times \ch{D}'_{s})$, where we have $\NL_{\tau,1,s} \simeq \ch{E}_{s} \times \ch{D}'_{s}$ with $\ch{D}'_{s}$ a factor isomorphic to $\ch{D}_{s}$ under $\rho$. (We recall that $\NL_{\tau,1,s}$ is a flag variety that factors in accordance with the factorization $\mathbf{G}_{\tau,s} = \mathbf{N}_{S,s} \cdot \mathbf{H}_{S,s}$.) Since all such varieties are parameterized algebraically by a Noetherian parameter space whose points correspond to the possible choices for the tuple $(s, \iota, L^{\bullet}_{c}, q)$, we can construct the desired $g$. After replacing $\mathcal{M}$ with a locally closed subscheme, we can assume that the fibres of $g$ all have dimension $a$.

Now apply \autoref{zariskidensityprop} to the family $g$ with $d = e+1$, and replace $R$ with $R[1/r!]$ (with $r$ as in the statement of \autoref{zariskidensityprop}). We show that the points of the loci $[S_{\kappa(\mathfrak{p}),\alpha}]^{>e}$ lie in the closed subscheme $E \subset S$ thus obtained for each $\mathfrak{p} = \mathfrak{p}_{j}$. Letting $\kappa = \kappa(\mathfrak{p})$ it suffices to show that $C := ([S_{\overline{\kappa},\alpha}]^{>e})^{\textrm{red}} \subset E_{\overline{\kappa}}$, where we use the notation $(-)^{\textrm{red}}$ to denote the underlying (geometrically) reduced locus. Choosing an arbitrary point $s \in S(\overline{\kappa})$ lying in the smooth locus of some component of $C$, it then suffices to show that $s$ lies in $E(\overline{\kappa})$. Since the corresponding component of $C$ is formally smooth of dimension at least $d$ at $s$, the quotient $\mathcal{O}_{C, s} / \mathfrak{m}^{r+1}$, where $\mathfrak{m}$ is the maximal ideal corresponding to $s$, admits a non-degenerate jet $j : \mathbb{D}^{d}_{r,\overline{\kappa}} \to S_{\overline{\kappa}}$ which factors through $\Spec \mathcal{O}_{C, s}$. Let $(\psi, \iota)$ be a framed $r$-limp for the data $\restr{(\mathcal{H}, F^{\bullet}, \nabla)}{S_{\overline{\kappa}}}$ which is defined at $s$. From \autoref{factorthroughlem} proven above, the composition $\psi \circ j$ lands inside $\ch{\iota}(\overline{\mathcal{S}(\mathfrak{p},\alpha)}^{\textrm{Zar}})$, which is a fibre of the map $g$. Thus, the data $(\iota, j, \psi \circ j)$ defines a point of $\mathcal{E}_{r}$ (in the notation of \autoref{zariskidensityprop}), and therefore the point $s$, which is the base-point of $j$, lies in $E = \overline{\textrm{im}(\mathcal{E}_{r} \to S)}^{\textrm{Zar}}$. 

Now let $F$ be the Zariski closure in $S$ of $\bigcup_{j = 1}^{\infty} [S_{\kappa(\mathfrak{p}_{j}), \alpha}]^{>e}$, which we have shown is a proper Zariski closed subset of $S$. Let $F_{1} \subset F$ be a component mapping to the generic point of $\Spec R$. Note that if $F^{\textrm{sm}}_{1} \subset F_{1}$ is the smooth locus, then $(\bigcup_{j = 1}^{\infty} [S_{\kappa(\mathfrak{p}_{j}), \alpha}]^{>e}) \cap F^{\textrm{sm}}_{1}$ is Zariski dense in $F^{\textrm{sm}}_{1}$. If at least one component of $F_{1,\mathbb{C}}$ does not lie inside a proper weakly special subvariety, then the algebraic monodromy group on any geometric component of $F^{\textrm{sm}}_{1,\mathbb{C}}$ agrees with the algebraic monodromy of $\mathbf{H}_{S}$ (here we use the fact that a conjugate of a weakly special subvariety is weakly special, see \cite[Theorem 2.6(b)]{fieldsofdef}). Thus, after restricting to the family $f : X_{F^{\textrm{sm}}_{1}} \to F^{\textrm{sm}}_{1}$ obtained by base-change we would be able to replace $S$ with $F^{\textrm{sm}}_{1}$ and repeat the above argument on some component of $F^{\textrm{sm}}_{1}$ to arrive at a contradiction (a smaller Zariski closed locus than $F$ containing the desired union); here we note that both $d$ and $e$ would be unchanged upon restriction, as would be $\ch{D}_{\mathbb{C}}$. We conclude that each component of $F_{\mathbb{C}}$ lies in a weakly special subvariety of $S_{\mathbb{C}}$ for the variation $\mathbb{V}$. 
\end{proof}

\section{Application in the Shimura Setting}
\label{shimappsec}

We now prove \autoref{genshimcase}, which implies \autoref{abvarcase} formally. Note that thus far we have not actually defined Bruhat strata for (special fibres of integral canonical models of) Shimura varieties, so our first task is actually to make sense of the the statement of \autoref{genshimcase}. Definitions of such strata already exist in the literature for at least the PEL type case, for instance in \cite{bruhatPELtype}, and our definition will be compatible with these existing definitions. Moreover, there exist already definitions of the Ekedahl-Oort strata of special fibres of integral canonical models of abelian type, as for instance given in \cite{shen2021stratifications}, and from this and the relation between Ekedahl-Oort and Bruhat strata in the PEL case one can infer a certain ``correct'' definition, which is the one we will give. We claim no originality in this definition, which is essentially a modification of the definitions appearing in \cite{shen2021stratifications}, but we do not know if it has thus far been explicitly stated. 

We begin with the case of Shimura varieties of Hodge type, mirroring the discussion in \cite{shen2021stratifications}. Let $(G, D)$ be a Shimura datum of Hodge type for which $G_{\mathbb{Q}_{p}}$ extends to a reductive group scheme $G_{\mathbb{Z}_{p}}$ over $\mathbb{Z}_{p}$. It will suffice for our purposes to explain the construction of integral canonical models for abelian type Shimura varieties $\textrm{Sh}_{K}(G,D)$ with $K = K_{p}K^{p}$ a neat arithmetic subgroup with $K_{p} \subset G(\mathbb{Q}_{p})$ hyperspecial, and $K^{p} \subset G(\mathbb{A}^{p}_{f})$. As explained in \cite[\S1.1]{shen2021stratifications} we can arrange for all of the following:
\begin{itemize}
\item[-] we have a symplectic embedding $i : (G, D) \hookrightarrow (\textrm{GSp}(V, \psi), D')$, where $V$ is a $\mathbb{Q}$-vector space of dimension $2g$, and $\psi$ is a non-degenerate symplectic form;
\item[-] $i_{\mathbb{Q}_{p}} : G_{\mathbb{Q}_{p}} \hookrightarrow \GL(V_{\mathbb{Q}_{p}})$ extends to a unique closed embedding $G_{\mathbb{Z}_{p}} \hookrightarrow \GL(V_{\mathbb{Z}_{p}})$, and so in particular there is a $\mathbb{Z}$-lattice $V_{\mathbb{Z}} \subset V$ such that $G_{\mathbb{Z}_{(p)}}$, the Zariski closure of $G$ in $\GL(V_{\mathbb{Z}_{(p)}})$, is reductive;
\item[-] there is a closed embedding 
\[ e : \textrm{Sh}_{K}(G,D) \hookrightarrow \textrm{Sh}_{K'}(\textrm{GSp}(V, \psi), D')_{E} , \]
with $E = E(G,X)$ the reflex field of $(G,D)$, and $K' \subset \textrm{GSp}(V,\psi)(\mathbb{A}_{f})$ a subgroup such that $\textrm{Sh}_{K'}(\textrm{GSp}(V, \psi), D')$ is a moduli scheme $\mathcal{A}_{g,1,K'}$ for abelian varieties with level-$K'$ structure;
\item[-] the integral canonical model $\mathcal{X} = \mathcal{I}_{K}(G,D)$ of $\textrm{Sh}_{K}(G,D)$ is the Zariski closure of the image of $e$ inside $\mathcal{A}_{g,1,K'} \times \Spec \mathcal{O}_{E,(v)}$, where $v$ is a fixed place of $E$ lying above $p$ (the usual normalization step is not required in light of \cite{normnotreq}).
\end{itemize}

Consider the universal abelian scheme $\mathcal{B}' \to \mathcal{A}_{g,1,K}$, and its pullback to an abelian scheme $\pi : \mathcal{B} \to \mathcal{X}$. Associated to this abelian scheme we obtain a triple $(\mathbb{V}, F^{\bullet}, \mathcal{H})$ with $\mathbb{V} = R^{1} \pi_{*} \mathbb{Z}$, $F^{\bullet}$ the Hodge filtration, and $\mathcal{H} = R^{1} \pi_{*} \Omega_{\mathcal{B}/\mathcal{X}}$, as in \S\ref{hodgetheoreticintro}. The special fibre of $\mathcal{H}$ at $v$ also comes with a natural conjugate filtration, again following \S\ref{hodgetheoreticintro}. As explained in \cite[\S1.1]{shen2021stratifications}, there are natural tensors $\nu_{1}, \hdots, \nu_{b}$ which give global, flat sections of tensor powers of $\mathcal{H}$ and its dual, and are identified with global Hodge tensors of $\mathbb{V}$ under the de Rham-Betti comparison. These tensors serve to define realizations of $G$ in each fibre of $\mathcal{H}$. We moreover have that:

\begin{lem}
Each tensor $\nu_{i}$ satisfies the property ($\mathfrak{p}_{v}$-Hodge) of \S\ref{hodgetheoreticintro}, where $\mathfrak{p}_{v} \subset \mathcal{O}_{E,(v)}$ is the prime corresponding to $v$. 
\end{lem}

\begin{proof}
Write $\mathcal{H}_{\textrm{cris}} = R^{1} f_{*,\textrm{cris}} \mathcal{O}_{\mathcal{B}/W(\kappa)}$ for the natural crystal over the crystalline site of the special fibre $\mathcal{X}_{\kappa}$, where $\kappa$ is the residue field of $\mathcal{O}_{E,(v)}$. The result will follow from \cite[Cor. 6.5]{algcycleloci} if we can show that each $\nu_{i}$ extends to a global tensor of $\bigoplus_{a,b \geq 0} \mathcal{H}^{\otimes a}_{\textrm{cris}} \otimes (\mathcal{H}^{*}_{\textrm{cris}})^{\otimes b}$ which is invariant under the geometric Frobenius endomorphism in each fibre. In the language of Dieudonn\'e crystals, this is \cite[Prop. 0.2]{kisin2010integral}.
\end{proof}

\begin{defn}
The Bruhat stratification of $\mathcal{X}_{\kappa}$ is the stratification of $\mathcal{X}$ induced by the stratification of $\mathcal{X}(\overline{\kappa})$ by the isomorphism class of the quadruple $(\mathcal{H}_{s}, F^{\bullet}_{s}, F^{\bullet}_{c,s}, \{ \nu_{i,s} \}^{b}_{i=1} )$, where $s \in \mathcal{X}(\overline{\kappa})$ is a point. 
\end{defn}

\noindent One may check that this definition is independent of the choices made in the construction (e.g., the symplectic embedding $i : (G, D) \hookrightarrow (\textrm{GSp}(V, \psi), D')$) as is done for instance in \cite{EOindep} in the case of Ekedahl-Oort strata; we omit the details.

\vspace{0.5em}

Now let us turn to the case of integral canonical models for Shimura varieties of abelian type, following \cite[\S1.2]{shen2021stratifications}. We will actually be able to omit most of the details of the actual construction, as we will see that we will need only the formal relationship between the Hodge and abelian-type cases to prove our main result. 

Recall that if $(G, D)$ is an abelian type Shimura datum, then there exists a Hodge type Shimura datum $(G_{1}, D_{1})$ such that one has a central isogeny $G^{\textrm{der}}_{1} \to G^{\textrm{der}}$ which induces an isomorphism of Shimura data $(G^{\textrm{ad}}_{1}, D^{\textrm{ad}}_{1}) \xrightarrow{\sim} (G^{\textrm{ad}}, D^{\textrm{ad}})$. Moreover, as explained in \cite[\S1.2.5]{shen2021stratifications} if $G_{\mathbb{Q}_{p}}$ extends to a reductive group $G_{\mathbb{Z}_{p}}$ over $\mathbb{Z}_{p}$, one can arrange for the existence of a reductive model $G_{1,\mathbb{Z}_{p}}$ of $G_{1}$, and even find models $G_{\mathbb{Z}_{(p)}}$ and $G_{1,\mathbb{Z}_{(p)}}$ such that we have an isomorphism $G^{\textrm{der}}_{1,\mathbb{Z}_{(p)}} \xrightarrow{\sim} G^{\textrm{der}}_{\mathbb{Z}_{(p)}}$ extending the one at the generic fibre. 

Letting now $K_{p} = G_{\mathbb{Z}_{p}}(\mathbb{Z}_{p})$ and $K_{1,p} = G_{1,\mathbb{Z}_{p}}(\mathbb{Z}_{p})$, let $\textrm{Sh}_{K_{p}}(G,D)^{+} \subset \textrm{Sh}_{K_{p}}(G,D)$ (resp. $\textrm{Sh}_{K_{1,p}}(G_{1},D_{1})^{+} \subset \textrm{Sh}_{K_{1,p}}(G_{1},D_{1})$) be the connected component corresponding to $1 \in G(\mathbb{A}_{f})$ (resp. $1 \in G_{1}(\mathbb{A}_{f})$). Following Deligne \cite{delignecanon} there exists a finite group $\Delta$ acting freely on $\textrm{Sh}_{K_{1,p}}(G_{1},D_{1})^{+}$ and an identification $\textrm{Sh}_{K_{1,p}}(G_{1},D_{1})^{+} / \Delta = \textrm{Sh}_{K_{p}}(G,D)^{+}$. In Kisin's integral model theory \cite{kisin2010integral} (c.f. \cite[\S1.2.5]{shen2021stratifications} this description, in addition to the group $\Delta$, is extended to the integral level over the ring $O_{F,(p)}$; here $F$ is the maximum unramified extension of the reflex field $E_{1}$ of $(G_{1}, D_{1})$, and $O_{F,(p)}$ is the localization at $p$ of its ring of integers. In particular, one can construct the integral canonical model of $\textrm{Sh}_{K_{p}}(G,D)^{+}$ as $\mathcal{I}_{K_{p}}(G,D)^{+} = \mathcal{I}_{K_{1,p}}(G_{1},D_{1})^{+} / \Delta$. 

With this description in mind, one can show the Bruhat stratification of $\mathcal{I}_{K_{1,p}}(G_{1},D_{1})^{+}$, which is induced by the Bruhat stratification of the integral schemes at each level, is invariant under the free action by $\Delta$. This is done in \cite[Prop. 3.4.3]{shen2021stratifications} in the setting of Ekedahl-Oort strata, which also implies the same statement for Bruhat strata once one shows that each Bruhat stratum is a finite union of Ekedahl-Oort strata. This latter fact, in turn, is easily seen from the definition of Ekedahl-Oort strata in terms of $F$-zips with $G$-structure given in \cite[\S3]{shen2021stratifications}, since the Ekedahl-Oort strata classify the same data as the Bruhat strata except with the inclusion of an additional isomorphism $\textrm{gr} F^{\bullet} \xrightarrow{\sim} \textrm{gr} F^{\bullet}_{c}$ of associated graded filtrations induced by Frobenius. One therefore finally defines the Bruhat strata of $\mathcal{I}_{K_{p}}(G,X)^{+}$, and hence of each constituent scheme in the associated tower, as the image of the Bruhat strata under the quotient map $\mathcal{I}_{K_{1,p}}(G_{1},D_{1})^{+} \to \mathcal{I}_{K_{1,p}}(G_{1},D_{1})^{+} / \Delta$.

\begin{proof}[Proof of \ref{genshimcase}]
Let $(G,D)$ be the associated abelian type Shimura data. Let us note that although the above constructions are carried out over localized rings and fields after fixing a prime $p$ and a place $v$ of the reflex field $E = E(G,D)$, such constructions can also be carried out over number rings $\mathcal{O}_{E}[1/N]$ for an integer $N$ for which $G$ admits a reductive model $\mathcal{G}$ over $\mathbb{Z}[1/N]$; this is checked explicitly in \cite{lovering2017integral}. (Note that such a model for $G$ can always be found by spreading out.) In particular, one can always find a Hodge-type Shimura datum $(G_{1}, D_{1})$ mapping to $(G, D)$ and an integral model $\mathcal{X}_{1}$ such that a finite map $X_{1} \to X$ of associated Shimura varieties spreads out to a finite map $\mathcal{X}_{1} \to \mathcal{X}$, agreeing with the canonical map constructed at all but finitely many primes. Because special subvarieties of $X_{1,\mathbb{C}}$ map to special subvarieties of $X_{\mathbb{C}}$ and Bruhat strata map to Bruhat strata under this map, this formally reduces the statement of \autoref{genshimcase} from the abelian type case to the Hodge type case.

In the case where $(G, D)$ has Hodge type we have already seen that Bruhat strata are defined exactly as in \S\ref{hodgetheoreticintro} with respect to cohomological data $(\mathcal{H}, F^{\bullet}, \{ \nu_{i} \}_{i=1}^{b})$ defined on $\mathcal{X}$ and the conjugate filtration at each prime. The result will then follow from \autoref{bruhatcritthm} if we can show that $\dim \ch{D}_{\overline{\kappa(\mathfrak{p}_{j})}, \alpha} = \dim \mathcal{X}_{\alpha, \mathfrak{p}_{j}}$. We have seen in \autoref{bruhatdim} that $\dim \ch{D}_{\overline{\kappa(\mathfrak{p}_{j})}, \alpha}$ is given by the length of the longest Weyl group element in the double coset $W_{o} w_{\alpha} W_{1}$, with the notation as in \S\ref{stratsec}. The result will then follow if we can show that this is also the dimension of the largest Ekedahl-Oort stratum contained inside $\dim \mathcal{X}_{\alpha, \mathfrak{p}_{j}}$. But the Ekedahl-Oort strata contained in the fixed Bruhat stratum corresponding to $W_{o} w_{\alpha} W_{1}$ are indexed by the cosets $W_{o} w_{\alpha} w_{1}$ as $w_{1}$ ranges over $W_{1}$ (see \cite[\S3.2]{shen2021stratifications}). The result therefore follows from \cite[Theorem B (ii)]{shen2021stratifications}, which shows that Ekedahl-Oort strata have dimension equal to the length of the longest vector in their corresponding coset. 
\end{proof}

\bibliography{hodge_theory}

\begin{thebibliography}{CPMS03}

\bibitem[Bor91]{borel1991linear}
A.~Borel.
\newblock {\em Linear Algebraic Groups}.
\newblock Graduate Texts in Mathematics. Springer New York, 1991.

\bibitem[BT17]{AXSCHAN}
Benjamin Bakker and Jacob Tsimerman.
\newblock The {A}x-{S}chanuel conjecture for variations of {H}odge structures.
\newblock {\em Inventiones mathematicae}, 217:77--94, 2017.

\bibitem[CPMS03]{CMS}
James~A Carlson, C.~(Chris) Peters, and Stefan M\"uller-Stach.
\newblock {\em Period mappings and period domains}.
\newblock Cambridge, U.K. : Cambridge University Press, 2003.

\bibitem[Del79]{delignecanon}
Pierre Deligne.
\newblock Vari\'et\'es de shimura: interpr\'etation modulaire, et techniques de
  construction de mod\'eles canoniques.
\newblock In A.~Borel and Casselman W., editors, {\em Automorphic forms,
  representations and L-functions (Corvallis 1977)}, pages 247--289. Amer.
  Math. Soc., 1979.

\bibitem[GGK12]{GGK}
Mark Green, Phillip Griffiths, and Matt Kerr.
\newblock {\em Mumford-Tate Groups and Domains: Their Geometry and Arithmetic
  (AM-183)}.
\newblock Princeton University Press, 2012.

\bibitem[Kat70]{katznil}
Nicholas Katz.
\newblock Nilpotent connections and the monodromy theorem : applications of a
  result of {Turrittin}.
\newblock {\em Publications Math\'ematiques de l'IH\'ES}, 39:175--232, 1970.

\bibitem[Kat72]{Katz1972}
Nicholas Katz.
\newblock Algebraic {S}olutions of {D}ifferential equations (p-curvature and
  the {H}odge {F}iltration).
\newblock {\em Inventiones mathematicae}, 18:1--118, 1972.

\bibitem[Kis10]{kisin2010integral}
Mark Kisin.
\newblock Integral models for shimura varieties of abelian type.
\newblock {\em Journal of the American Mathematical Society}, 23(4):967--1012,
  2010.

\bibitem[KO68]{katz1968}
Nicholas~M. Katz and Tadao Oda.
\newblock On the differentiation of {D}e {R}ham cohomology classes with respect
  to parameters.
\newblock {\em J. Math. Kyoto Univ.}, 8(2):199--213, 1968.

\bibitem[KO21]{2021InMat.225..857K}
B.~{Klingler} and A.~{Otwinowska}.
\newblock {On the closure of the {H}odge locus of positive period dimension}.
\newblock {\em Inventiones Mathematicae}, 225(3):857--883, September 2021.

\bibitem[KOU20]{fieldsofdef}
Bruno {Klingler}, Anna {Otwinowska}, and David {Urbanik}.
\newblock {On the fields of definition of Hodge loci}.
\newblock {\em To appear in Ann. Sci. de l'École Nor. Sup. arXiv:2010.03359},
  October 2020.

\bibitem[Lov17]{lovering2017integral}
Tom Lovering.
\newblock Integral canonical models for automorphic vector bundles of abelian
  type.
\newblock {\em Algebra \& Number Theory}, 11(8):1837--1890, 2017.

\bibitem[MPT19]{axschanshimura}
Ngaiming Mok, Jonathan Pila, and Jacob Tsimerman.
\newblock {A}x-{S}chanuel for {S}himura varieties.
\newblock {\em Annals of Mathematics}, 189(3):945--978, 2019.

\bibitem[SZ85]{Steenbrink1985}
Joseph Steenbrink and Steven Zucker.
\newblock Variation of mixed {H}odge structure. {I}.
\newblock {\em Inventiones mathematicae}, 80:489--542, 1985.

\bibitem[SZ21]{shen2021stratifications}
Xu~Shen and Chao Zhang.
\newblock Stratifications in good reductions of {S}himura varieties of abelian
  type, 2021.

\bibitem[Urb21]{periodimages}
David Urbanik.
\newblock On the {T}ranscendence of {P}eriod {I}mages.
\newblock {\em arXiv preprint arXiv:2106.09342}, 2021.

\bibitem[Urb22]{algcycleloci}
David Urbanik.
\newblock Algebraic {C}ycle {L}oci at the {I}ntegral {L}evel, 2022.

\bibitem[Wed14a]{bruhatPELtype}
Torsten Wedhorn.
\newblock {\em Mathematische Zeitschrift}, 277:725--738, 2014.

\bibitem[Wed14b]{bruhatstack}
Torsten Wedhorn.
\newblock {\em M\"unster Journal of Mathematics}, 7:529--556, 2014.

\bibitem[Xu20]{normnotreq}
Yujie Xu.
\newblock Normalization in integral models of shimura varieties of hodge type,
  2020.

\bibitem[Yve92]{Andre1992}
Andr\'e Yves.
\newblock Mumford-{T}ate groups of mixed {H}odge structures and the theorem of
  the fixed part.
\newblock {\em Compositio Mathematica}, 82(1):1--24, 1992.

\bibitem[Zha14]{EOindep}
Chao Zhang.
\newblock Remarks on {E}kedahl-{O}ort stratifications, 2014.

\end{thebibliography}
\bibliographystyle{alpha}

\end{document}